\documentclass[conference]{IEEEtran}
\usepackage{cleveref}
\usepackage{ifpdf}

\usepackage[noadjust]{cite}

\ifCLASSINFOpdf
  \usepackage[pdftex]{graphicx}
   \graphicspath{{../images/},{../../SIMU/plots/keep/}}
   \DeclareGraphicsExtensions{.pdf,.jpeg,.png,.svg}
\else
   \usepackage[dvips]{graphicx}
   \graphicspath{{../eps/}}
   \DeclareGraphicsExtensions{.eps}
\fi
\usepackage{amsmath}
\usepackage{frenchineq}

\usepackage{algorithmic}
\usepackage{bm}
\usepackage{textcomp}
\usepackage{subfig}

\newcommand{\dssum}{\displaystyle\sum}
\newcommand{\SC}{\text{S\hspace{-1pt}C}}

\newcommand{\eqDef}{\overset{\text{def}}{=}}

\newtheorem{theorem}{Theorem}

\newtheorem{definition}{Definition}

\Crefname{equation}{Eq.}{Eqs.}
\Crefname{figure}{Fig.}{Figs.}
\Crefname{tabular}{Tab.}{Tabs.}
\Crefname{table}{Tab.}{Tabs.}
\Crefname{theorem}{Thm.}{Thms.}
\Crefname{definition}{Def.}{Defs.}

\usepackage[colorinlistoftodos,bordercolor=orange,backgroundcolor=orange!20,linecolor=orange,textsize=scriptsize]{todonotes}
\usepackage{multirow}
\usepackage{pbox}

\begin{document}

\title{Demand Side Management in the Smart Grid: \\ an Efficiency and Fairness Tradeoff}



%
\author{\IEEEauthorblockN{Paulin Jacquot\IEEEauthorrefmark{1}\IEEEauthorrefmark{2} \IEEEmembership{Student Member, IEEE},
Olivier Beaude\IEEEauthorrefmark{1}, 
St\'ephane Gaubert\IEEEauthorrefmark{2}, Nadia Oudjane\IEEEauthorrefmark{1}}
\IEEEauthorblockA{\IEEEauthorrefmark{1} EDF, Saclay, France\\
}
\IEEEauthorblockA{\IEEEauthorrefmark{2} Inria Saclay, CMAP, Ecole Polytechnique, Palaiseau, France}  \thanks{Accepted for presentation at IEEE International Conference on Innovative Smart Grid Technologies Europe (ISGT), 2017, Torino, Italy.}
}

\IEEEoverridecommandlockouts


\maketitle

\begin{abstract}

We compare two Demand Side Management (DSM) mechanisms, introduced respectively by Mohsenian-Rad \textit{et al } (2010) and Baharlouei \textit{et al} (2012), in terms of efficiency 
 and fairness. 
Each mechanism defines a game where the consumers optimize their flexible consumption to reduce their electricity bills. Mohsenian-Rad \textit{et al} propose a daily mechanism for which they prove the social optimality. Baharlouei \textit{et al} propose a hourly billing  mechanism for which we give theoretical results:  we prove the uniqueness of an equilibrium in the associated game and give an upper bound on its price of anarchy. We evaluate numerically the two mechanisms, using real consumption data from \textit{Pecan Street Inc.} The simulations show that the equilibrium reached with the hourly mechanism is socially optimal up to $0.1\%$, and that it achieves an important fairness property according to a quantitative indicator we define. 
We observe that the two DSM mechanisms avoid the synchronization effect induced by non-game theoretic mechanisms, e.g. Peak/OffPeak hours contracts. 

\end{abstract}


%
\IEEEpeerreviewmaketitle

\section{Introduction}

Demand Side Management (DSM) is a way to exploit demand elasticity to achieve energy balancing or ancillary services \cite{saad2012game}. It has become an increasing issue for electrical systems with the new possibilities offered by smart grid technologies. 
Such flexibilities  can be an efficient way 
 of reducing peak electricity demand. 
These opportunities could enable to avoid using expensive flexible power plants, with high variational costs and important air emissions.  
Besides, with an increasing part of fluctuating renewable production sources, 
flexibilities 
will be crucial to ensure the network reliability. Thus, DSM could provide significant economic, reliability and environmental benefits.

Demand Response has been the subject of a blooming literature in the past decade. One can refer to \cite{siano2014demand}, \cite{vardakas2015survey} and \cite{saad2012game} for surveys on this topic. Several works \cite{mohsenian2010autonomous}, \cite{chen2010two}, \cite{ibars2010distributed},\cite{li2011optimal}  propose game theory frameworks and dynamic pricing models that enable  users to schedule their flexible electricity consumption.

In this paper, we consider the framework proposed by Mohsenian-Rad \textit{et al} (2010) \cite{mohsenian2010autonomous}. 
 An isolated 
operator holding power production assets 
 proposes to equip each of its consumers 
 with an automatic Energy Consumption Scheduler (ECS), integrated inside a smart meter. 
This ECS is connected both to the power grid and to flexible electrical appliances such as plug-in electric vehicles, air conditioning, heating, etc. 
Every day, each ECS  runs locally an algorithm 
to schedule the local consumption for each hour of the next day. This scheduling aims to minimize the total energy cost for the system while respecting the quality of service and constraints for each flexible appliance. Each consumer pays a part of the total generation cost  proportional to the energy it has consumed in the day. 
Technically, the distributed algorithm consists in implementing a Best Response Dynamic: each ECS iteratively minimizes its individual bill, parametrized by its information on the total load of the system. Then, it communicates to the system its individual optimal load. The system updates the current total load and informs all ECSs. 

The implementation of a DSM mechanism raises several difficulties.
Owing to the huge number of variables and constraints and to the impossibility for an aggregator to collect all the consumption constraints because of privacy concerns, the optimization has to rely on a decentralized algorithm that minimizes the information exchanged with the users. 
%
%
Obviously, efficiency is also a desirable property. One wishes that, as in \cite{mohsenian2010autonomous}, the scheduling process leads to an optimal or close to optimal consumption profile for the global system costs, while respecting all the users constraints. 
%
Another important feature, more discussed in \cite{baharlouei2012tackling}, is fairness:  
the payment model should penalize consumers imposing costly constraints for the system, while rewarding flexible ones. 
This point is essential to ensure the acceptability of the process and encourage users to stay in a DSM program. Besides, this feature has a merit in terms of incentives. Indeed, such fair billing models would encourage consumers to modify their constraints so that available flexibility for the system increases.
%

The contributions of this paper are twofold. First, we give new theoretical results associated with the hourly billing mechanism \cite{baharlouei2012tackling} by proving the uniqueness of a Nash Equilibrium, and by specifying an explicit upper bound on its price of anarchy. Next, we present numerical results, based on simulations using real consumption data, that compare the two billing mechanisms \cite{mohsenian2010autonomous} and \cite{baharlouei2012tackling} in terms of efficiency and fairness. The results show that the hourly mechanism achieves a very small price of anarchy and an important fairness property.

The paper is organized as follows. In \Cref{sec: energy cons game} we introduce the consumption game model. In \Cref{sec: indicators} we recall some notions from game theory 
and define quantitative indicators to measure the efficiency and the fairness of a given billing mechanism. In \Cref{sec:DP} we introduce the daily proportional billing considered in \cite{mohsenian2010autonomous} and recall its main properties.  In \Cref{sec:HP}, we focus on the hourly proportional billing proposed in \cite{baharlouei2012tackling} and present our new theoretical results. Last, \Cref{sec: simu} is devoted to numerical experiments, based on real consumption data, from which we compare the fairness and efficiency of the two different billing mechanisms.
%
%
%
%
%
%
%
%
%
%
%
%
%
%
%
%
%
%
%
%
%
%
%
%
\section{Energy Consumption Game}
\label{sec: energy cons game}

We consider an autonomous network composed of a unique electricity provider and a set $\mathcal{N}$ of $N$ electricity consumers. We use a model similar to Mohsenian-Rad \textit{et al} \cite{mohsenian2010autonomous}.

\subsection{Consumers constraints}

Each user $n$ has a set of electric appliances $\mathcal{A}_n$. For each $a\in\mathcal{A}_n$, this user (himself or through an ECS) can set the power  $x_{na}^h$  allowed to $a$ at each time period $h$ in $\mathcal{H}_{na}=\{\alpha_{na},\dots,\beta_{na}\} \subset \mathcal{H}$, where $\mathcal{H}$ is the set of time periods considered over a day. We consider the following constraints:
\vspace{-10pt}
\begin{subequations}
\label{eq:userConstraints}
\begin{align}
\label{cons: HP total power} & \hspace{0.5cm}   \dssum_{h\in \mathcal{H}} x^h_{na} = E_{na} , \ \forall a \in \mathcal{A}_n \ , \\ 
 \label{cons: HP minmax power}& \hspace{0.5cm} \underline{x}^h_{na} \leq x_{na}^h 
\leq \overline{x}^h_{na} , \ \forall a \in \mathcal{A}_n, \forall h \in \mathcal{H} \ .
  \end{align}
\end{subequations}
Each electric appliance $a\in \mathcal{A}_n$ requires a fixed daily amount of energy (\ref{cons: HP total power}). Due to physical limits, the power set to appliance $a$ is bounded from below and above (\ref{cons: HP minmax power}). If the appliance $a$ can not be used in the time period $h$, we set $\underline{x}^h_{na} = \overline{x}^h_{na} =0$. The set of available time periods for $a\in\mathcal{A}_n$ is therefore given by $\mathcal{H}_{na}=\{h : \ \overline{x}_{na}^h >0 \}$.

We will denote more compactly by $\mathcal{X}_n$ the set of feasible loads $(x_{na}^ h)_{h,a} $ that respect the constraints given by (\ref{eq:userConstraints}).

\subsection{System cost functions}

We denote by $C_h(\ell^h)$ the system costs for providing to users the total load $\ell^h:=\sum_{a,n} x_{na}^h$ at time $h$. It is widely accepted that marginal production costs increase with demand. Hence we assume that cost functions $C_h(.)$ are \emph{increasing} and \emph{strictly convex} \cite{mohsenian2010autonomous}. These functions can depend on the period $h$ as it is more expensive to produce energy on peak hours or for instance, when the renewable production is low. Practically, in our model, $(C_h)_{h\in\mathcal{H}}$ can be the actual costs for the provider but can also be an artificial signal that is sent to each user's ECS in order to make him perform a decentralized optimization of his consumption.
We will mostly consider quadratic cost functions as done in \cite{mohsenian2010autonomous}:
\begin{equation}
\label{eq:quad Ch}
C_h(\ell^h)=a_2^h {(\ell^h)}^2+ a_1^h \ell^h + a_0^h \ .
\end{equation}
We assume for simplicity that the total system cost, denoted by $\mathcal{C}:=\sum_h C_h$ is divided among consumers, in a way defined by the provider.  If we denote by $b_n$ the bill paid by user $n$ for the day, then we have $\mathcal{C}=\sum_n b_n$. In general, $b_n$ can depend on the induced costs $(C_h)_h$ and the load vector of each user $(\ell_n^h)_{n,h}$, where $\ell_n^h:=\sum_{a}x^h_{na}$. In this work we are interested in two different billing models given in Sections \ref{sec:DP} and \ref{sec:HP}. 

Finally, each user will try to minimize his bill $b_n$ while respecting his constraints \eqref{eq:userConstraints}, by solving the problem:
\begin{equation}
\label{eq:userMinimize}
\min_{\bm{x}_n \in \mathcal{X}_n} b_n(\bm{x}_n,\bm{x}_{-n})
\end{equation}
where $\bm{x}_n := (x_{na}^h)_{h,a} $ and $\bm{x}_{-n}=(x_m)_{m \neq n}$ stands for the consumption vector of all users but $n$. As $b_n$  depends on both $\bm{x}_n$ and $\bm{x}_{-n}$, this problem can be stated in the framework of game theory. We refer the reader to \cite{fudenberg1991game} for background. With $\mathcal{N}$ denoting the set of players, $\mathcal{X}:=\prod_{n\in \mathcal{N}} \mathcal{X}_n$ the set of pure strategies and $(b_n)_n$ the vector of bills, the game is formulated under normal form $ \mathcal{G}=\left(\mathcal{N},\mathcal{X},(b_n)_n\right)$.
%

\section{Measuring efficiency and fairness}
\label{sec: indicators}
\subsection{Efficiency and the Price of Anarchy}
\label{sec:Efficiency}

We define the social cost of a load solution $\bm{x}=(\bm{x}_n)_{n\in_\mathcal{N}}$ as the sum of the bills in the population, that is:
\begin{equation}
\SC( \bm{x}) := \sum_{n\in\mathcal{N}} b_n(\bm{x}_n,\bm{x}_{-n}) \ .
\end{equation}
Since we assume that the total system costs $\mathcal{C}$ are shared among the users, we also have the equality:
\begin{equation}
\mathcal{C}(\bm{\ell}):=\sum_{h\in\mathcal{H}} C_h(\ell^h)= \SC( \bm{x}) \ ,
\end{equation}
with $\ell^h:=\sum_{a,n}x^h_{na}$  the total load at time $h$ and $\bm{\ell}:=(\ell^h)_h$.

The efficiency of a mechanism is usually measured in game theory by the ratio of the optimum social cost of the system $\SC^* := \inf_{\bm{x}\in\mathcal{X}} \SC(\bm{x})$ and the social cost of the worst Nash Equilibrium: 

\begin{definition}[\cite{koutsoupias1999worst}] {\emph{Price of Anarchy.}}
\label{def:PoA}

Given a game $ \mathcal{G}$ and $\mathcal{X}_{\mathcal{G}}^{\text{NE}}$ its set of Nash Equilibria
, the price of anarchy of $\mathcal{G}$ is given as:
\begin{equation}
\label{eq:PoA}
\text{PoA}(\mathcal{G}):=\frac{ \sup_{\bm{x}\in \mathcal{X}_{\mathcal{G}}^{\text{NE}}} \SC\left(\bm{x}\right)}{\SC^*} \ .
\end{equation}
\end{definition}

The notion of Price of Anarchy has been widely studied in congestion and routing games (see \cite{roughgarden2016twenty},\cite{johari2005efficiency},\cite{johari2006scalable}). Theoretical bounds have been established in particular frameworks (e.g. congestion games in \cite{christodoulou2005price},\cite{roughgarden2006potential},\cite{roughgarden2015intrinsic}).

\subsection{Fairness and Marginal Cost Pricing}
\label{sec:Fairness}

To design a fair mechanism, the bill $b_n$ paid by each user $n$ should reflect the cost user $n$ induces to the system, what we call the externality of $n$. Precisely, we denote by:
\begin{equation}
\mathcal{C}^*_{\mathcal{M}} := \displaystyle\inf_{(\bm{\ell}_m)_{m\in\mathcal{M}}} \textstyle{\sum}_{h\in\mathcal{H}} C_h\left( \sum_{m\in\mathcal{M}} \ell^h_m \right) 
\end{equation} 
the optimal system cost that can be achieved with the users in the set $\mathcal{M}$ while respecting their constraints.  The externality of user $n$ is the difference between the optimal system cost achieved with $n$ in the population and the optimal system cost that can be achieved without $n$, that is, $V_n:= \mathcal{C}^*_{\mathcal{N}}- \mathcal{C}^*_{\mathcal{N}\setminus \{n\}}$. The quantity $V_n$ is not necessarily proportional to the total energy the user asks per day, as the load distribution between peak and off-peak hours also impacts the system cost.

 To be fair, the bill of user $n$ should be proportional to $V_n$ \cite{baharlouei2012tackling}. This motivates the introduction of the following mechanism:
\begin{equation}
\label{eq: VCG billing}
b_n^{\textsc{vcg}}(\bm{x}_n,\bm{x}_{-n}):= \sum_{h\in\mathcal{H}} C_h\left(\sum_{m\in\mathcal{N}} \ell^h_m \right)- \mathcal{C}_{\mathcal{N}\setminus \{n\}}^*
\end{equation}
which, as noticed in \cite{samadi2012advanced}, corresponds to a Vickrey-Clarke-Groves mechanism (VCG, see \cite{clarke1971multipart}).  In particular, it minimizes the system cost, which implies that in this model at an equilibrium $\bm{x}^\text{NE}$ we will have $\forall n, b_n^{\textsc{vcg}}(\bm{x}^\text{NE})=V_n$. The authors in \cite{moulin2001strategyproof} defined, in a more general framework, this pricing as Marginal Cost Pricing, and showed (Prop. 3) that it is the unique VCG mechanism that satisfies reasonable conditions. However, the mechanism \eqref{eq: VCG billing} does not recover the system cost $\mathcal{C}^*_\mathcal{N}$, and should be renormalized as $b^\text{F}_n:= \frac{b^{\textsc{vcg}}_n}{\sum_m V_m} \mathcal{C}^*_{\mathcal{N}}$. Although being centralized and hardly tractable, the billing mechanism $b^\text{F}_n$ is efficient (PoA=1) and fair ($b^\text{F}_n \propto V_n$) and we take it as a reference, following \cite{baharlouei2013achieving},\cite{baharlouei2014efficiency}, to define a fairness measure of any billing mechanism: 

\begin{definition}[\cite{baharlouei2014efficiency}]{\emph{Fairness Index.}}
\label{def:fairness}

The fairness index of a billing mechanism $(b_n)_n$ is its maximal normalized distance to $(V_n)_n$ (or equivalently to $(b^\text{F}_n)_n$) at a Nash Equilibrium:
\vspace{-10pt}
\begin{equation}
\label{eq: fairness ind}
F := \sup_{\bm{x}\in \mathcal{X}_{\mathcal{G}}^{\text{NE}}} \left[ \sum_{n\in \mathcal{N}} \left\lvert \frac{V_n}{\sum_{m\in\mathcal{N}} V_m} -\frac{b_n\left(\bm{x} \right) }{\sum_{m\in\mathcal{N}} b_m\left(\bm{x} \right)  }\right\rvert \right].
\end{equation}
\end{definition}

In \cite{baharlouei2014efficiency}, the authors notice the link between $V_n$ and the notion of Shapley Value (\cite{shapley1953value}) defined for cooperative games. However, since the Shapley Value is given by a combinatorial formula involving all possible coalitions within $\mathcal{N}$, it becomes quickly untractable as the cardinality of $\mathcal{N}$ grows. It is therefore more appropriate to use $(V_n)_n$ as reference.

\section{Daily Proportional Billing: social optimality}
\label{sec:DP}

\subsection{Daily Proportional (DP) billing: definition}
In this section, we recall the standard billing mechanism of \cite{mohsenian2010autonomous}. Consumers share the total cost of the system proportionally to their total consumption over the day. More precisely, if we denote by $E_n=\sum_{a\in\mathcal{A}_n} E_{na}$ the total energy needed by $n$, the bill of this user is:\vspace{-0.15cm}
\begin{equation}
\label{eq: DP}
b_n^{\text{DP}}(\bm{x}_n,\bm{x}_{-n}) = \frac{E_n}{\sum_{m\in\mathcal{N}}E_m} \sum_{h\in\mathcal{H}} C_h( \ell^h) \ .
\end{equation}

\subsection{Properties}

As all users minimize $\SC(\bm{x})$ up to a constant factor, several properties follow (detailed proofs are in \cite{mohsenian2010autonomous}). 
 First, \cite[Thm.~1]{mohsenian2010autonomous} ensures that a Nash Equilibrium (NE) exists and, as cost functions $(C_h)_h$ are assumed strictly convex, it is unique in terms of aggregated load $\left(\ell^h\right)_h$. 
The NE minimizes the social cost $\SC$ (\cite[Thm.~2]{mohsenian2010autonomous}). In order to compute the NE in the game, the authors in \cite{mohsenian2010autonomous} consider the implementation of Best Response Dynamic, that can be defined as follows:
\begin{definition}{ \emph{Best Response Dynamic (BRD). }}
\label{def: BRD}

In BRD, at each iteration $k$, a user $n$ is randomly chosen, and solves his local optimization problem (\ref{eq:userMinimize}) with knowledge of the load of others  $\bm{\ell }_{-n}^{(k-1)}$, taken as a  parameter. The resulting load $\bm{\ell}^{*}_n$ is used to update   $\bm{\ell }_{n}^{(k)}=\bm{\ell}^{*}_n$ and $\bm{\ell }_{-n}^{(k)}=\bm{\ell }_{-n}^{(k-1)}$.
\end{definition}

In practice, we only need the aggregated load $\sum_{m\neq n} \ell^h_m$ to solve user $n$'s problem, thus privacy of the users is preserved.
Again, due to the proportionality of the users objectives, the BRD will converge to the NE (\cite[Thm.~3]{mohsenian2010autonomous}). Finally, \cite[Thm.~4]{mohsenian2010autonomous} ensures that no user $n$ can reduce his bill $b_n$ by giving wrong information about his load $(\ell_n^h)_h$ during the process.

\section{Hourly Proportional Billing: Fairness}
\label{sec:HP}

\subsection{Hourly Proportional (HP) billing: definition}

Here, the total cost is divided between consumers at each time period, according to the load they asked at this time period. Intuitively, this enables to bring to each user the real cost of its demand, in particular during peak hours. More formally, the bill of user $n$ is:
\vspace{-6pt}
\begin{equation}
\label{eq:HP}
b_n^{\text{HP}}(\bm{x}_n,\bm{x}_{-n}) = \sum_{h\in\mathcal{H}}\frac{\ell_n^h}{\sum_{m\in\mathcal{N}}\ell^h_m}  C_h( \ell^h) \ .
\end{equation}

This billing mechanism was already formulated in \cite{baharlouei2012tackling}, \cite{baharlouei2013achieving}, \cite{baharlouei2014efficiency}. In the latter, the authors show numerically that it is fairer (according to indicator $F$ \eqref{eq: fairness ind}) than the payment \eqref{eq: DP}.

\subsection{Properties}

With payments \eqref{eq:HP}, the game $\mathcal{G}=\left(\mathcal{N},\mathcal{X},(b_n^\text{HP})_n\right)$ has the following properties.

\begin{theorem}
\label{prop:exist}Let $c_h(\ell^h):=\frac{1}{\ell^h}C_h (\ell^h)$ be the per-unit price of electricity. If $c'_h \geq 0$, \textit{i.e.} prices are \textit{increasing} with global load, then a Nash Equilibirum exists. If, in addition:
\vspace{-4pt}
\begin{equation}
\label{e-impliesuniqueness}
\forall h, \frac{ ({\ell^h})^2 } {\sum_n ({\ell^h_n})^2}> \left(\frac{\ell^h c''_h(\ell^h)}{2c_h'(\ell^h) }\right)^2 
\end{equation}
then the Nash Equilibrium is unique.
\end{theorem}

Note that if the load is close to uniform, we have $ \frac{\ell_n^h}{\ell^h}\simeq \frac{1}{N}$ so  ${ {\ell^h}^2 }/{\sum {\ell^h_n}^2} \simeq N$, and Assumption~\eqref{e-impliesuniqueness} is satisfied
as soon as the network has enough users.
\Cref{prop:exist} is obtained as a consequence of the results of Rosen \cite{rosen1965existence} on the existence and uniqueness of Nash Equilibria for $N-$persons concave games.
The main technical difficulty is to show the strict diagonal convexity condition of Rosen. If we define the gradient of users costs $g(\bm{x}) \eqDef  (\nabla_nb_n(\bm{x}))_{n\in \mathcal{N}}$, Rosen's condition is equivalent to having the Jacobian $J(g)(\bm{x})$  definite positive for all $\bm{x}$ admissible. We prove it by using a linear algebra perturbation theorem.
The details of the proof will be given elsewhere.

In general, the Nash Equilibrium does not achieve social optimality. However,
the following result provides a bound on the Price of Anarchy (\Cref{def:PoA}).
%

\begin{theorem}
\label{th: PoA bound}
In the quadratic case (\ref{eq:quad Ch}) with no constant coefficient ($\forall h, a_0^h=0$),
the price of anarchy is bounded:
\vspace{-4pt}
\begin{equation}
\text{PoA} \leq 1+ \frac{3}{4}\sup_{h\in\mathcal{H}} \frac{1}{1+{a_1^h}/{( a_2^h \overline{\ell}^h} ) } \ .
\end{equation} 
\end{theorem}


The proof of this result relies 
on the property of local smoothness introduced by Roughgarden and Schoppmann.

\begin{definition} [\cite{roughgarden2015local}]{ \emph{Local smoothness.}}

A cost minimization game  $\mathcal{G}=(\mathcal{N},\mathcal{X},(b_n)_n) $ is locally $(\lambda,\mu)$-smooth with respect to $y$ iff for all feasible outcome $x$:
\vspace{-8pt}
\[ \ \sum_{n\in\mathcal{N}} b_n(x)+ \nabla_n b_n(x)^T (y_n-x_n) \leq \lambda \SC(y) +\mu \SC(x) \ .\]
\end{definition}
\vspace{-5pt}
The authors in \cite{roughgarden2015intrinsic} show that if $\mathcal{G}$ is locally $(\lambda,\mu)$-smooth with respect to an optimal outcome $y$, then the PoA is less or equal than $\frac{\lambda}{1-\mu}$.
In our case, with $r_h=\frac{a_1^h}{a_2^h \overline{l}^h}$, we prove that the game is locally $(\lambda,\mu)$-smooth with respect to any $y$ with:
\[\mu=\sup_h \frac{-1+\sqrt{1+(1+r_h)^2}}{(1+r_h)^2} \text{ and } \lambda=\sup_h \frac{(1+r_h\mu)^2+\mu}{4(1+r_h)\mu}.\]


\section{Numerical Experiments}
\label{sec: simu}

We compare numerically the billing mechanisms DP (\ref{eq: DP}) and  HP (\ref{eq:HP}), based on the two criterias of efficiency (Def.~\ref{def:PoA}) and fairness (Def.~\ref{def:fairness}). We extracted a set $\mathcal{N}$ of 30 users from the database \textit{PecanStreet Inc.} (\cite{PecanStreet }), which gathers hundreds of disaggregated residential consumption profiles in Texas, U.S.
We use hourly timesteps so that $\mathcal{H}=\{0,1,,...,23\}$. We consider that each day, just before midnight, the flexible consumption of each user for the next day is computed as an equilibrium strategy of the game, implementing a BRD algorithm as in \Cref{def: BRD}, using the Smart Grid for exchanging the information $\bm{\ell}_{-n}$. We run simulations day by day on the set of 30 days
\footnote{To start simulations with a working day, we dismissed January, $1^\text{st}$.} $\mathcal{D} := \{ 02/01/2016, \dots, 31/01/2016 \}$.

\subsection{Flexible Appliances: Electric Vehicles and Heating}
\label{subsec: flex appl}
We study a population of residential consumers owning electric vehicles (EV) and electrical heating systems (furnace).
EVs present an important flexibility (\cite{beaude2016reducing}) since an EV remains plugged in while it is parked, and a smart charging can be automated without constraints for the user.
Similarly, the initial consumption profile of a heating system can be modified  without strong impact on the comfort of the household.

In our simulations, we consider a first case where EVs are the only flexible appliances, accounting for $20.4\%$ of an average daily global energy of 1014kWh, and a second case where furnaces are also considered as flexible appliances, increasing the part of flexible energy to $25.8\%$.
The remaining of each user's consumption is nonflexible, which we denote by $(\hat{\ell}_{n,\text{NF}}^h)_h$ (the hat stands for \emph{observed} data values). 
\Cref{fig: load data and Peak/offPeak}(top) shows the repartition between flexible and nonflexible load on a typical day. The nonflexible load is  more important on some hours than others, so that even with hourly uniform system costs \eqref{eq: coefffs uni}, these hours will have bigger marginal costs.

The users constraints (\ref{eq:userConstraints}) are evaluated as follows: we consider two types of days: $\mathcal{D}_1$ for  weekdays (Monday to Friday) and $\mathcal{D}_2$ for weekend days (Saturday and Sunday). For each type $\mathcal{D}_k$ of day and each user $n$, we suppose that appliance $a$ can be used at $h$ if it exists a day of type $\mathcal{D}_k$ where $a$ was on at $h$. More precisely, for a day of type $\mathcal{D}_k$,
$\mathcal{H}_{na}= \bigcup_{d \in \mathcal{D}_k} \{h \ : \ x_a^{d,h} >0 \}.$
For simplicity, we took the min power $\underline{x}_{na}^h$ equal to 0 and the max power  $\overline{x}_{na}^h$  equal to the maximal nonnegative value found on the data set $\overline{x}_{na}^h= \displaystyle\max_{h,d \in \mathcal{H}\times \mathcal{D}}\hat{x}_{na}^{d,h}$, if $h\in \mathcal{H}_{na}$, and 0 otherwise. 

\subsection{System Costs}

We consider that the provider costs $(C_h)_h$ are functions of the total load $L^h := \ell^h_{\text{NF}}+ \ell^h$. 
More precisely, they are given as $C_h(\ell^h)=\tilde{C}(L^h)$ where  $\tilde{C}$  is given in dollar cents (\textcent) by:
\begin{equation}
\label{eq: coefffs uni}
\forall h \in \mathcal{H}, \tilde{C}_h(L^h)=\tilde{C}(L^h):= 0.1 + 8L^h + 0.04 {(L^h)}^2 \ . 
\end{equation}

The average hourly nonflexible load on all days in $\mathcal{D}$ is:
$$\langle \hat{\ell}_{\text{NF}}\rangle := \frac{1}{|\mathcal{D}|\times |\mathcal{H}|}\sum_{d,h \in \mathcal{D} \times \mathcal{H}} \left( \sum_{n\in\mathcal{N}} \hat{x}_{n,\text{NF}}^{d,h} \right)=31.3\text{kWh}.$$ The coefficients in \eqref{eq: coefffs uni} are chosen arbitrarily but such that the price ${C_h(\langle \hat{\ell}_{\text{NF}}\rangle )}/{\langle \hat{\ell}_{\text{NF}}\rangle}$ given by\eqref{eq: coefffs uni} match the price proposed by the distributor \textit{CoServ} \cite{Coserv} of 
$8.5$\textcent/kWh for base contracts.

We assume that the nonflexible load $(\hat{\ell}_{n,\text{NF}}^h)_h$ is billed in a separate process (for instance, according to a baseline contract as defined in the next subsection). We apply the proposed billing mechanisms DP \eqref{eq: DP} and HP \eqref{eq:HP} on the flexible part $\ell^h$ only. Although the system costs $(\tilde{C}^h)_h $ are hourly uniform, the variation in the nonflexible load $\ell^h_{\text{NF}}$ over the hours induces a variation on the  \emph{cost of flexible load} over the hours\footnote{In practice, the provider could rely on a forecast of $\ell^h_{\text{NF}}$ for the next day instead of the observed value $\hat{\ell}^h_{\text{NF}}$, to compute the functions $(C_h)_h$.}:
\begin{equation}
\begin{split}
\label{eq: cost flexible load}
\forall h \in \mathcal{H}, C_h(l^h) &:=\tilde{C}(\ell^h_{\text{NF}}+\ell^h) - \tilde{C}(\ell^h_{\text{NF}}) \\&= (8+0.08\ell^h_{\text{NF}})\ell^h + 0.04 {(\ell^h)}^2 \ . 
\end{split}
\end{equation}

\subsection{Two Reference Non-Game Theoretic Billing Models}
In order to compare the formulated game-theoretic models \eqref{eq: DP} and \eqref{eq:HP} to existent non game-theoretic billing models, we also consider the two following standard models as references:

1) \textbf{Baseline billing.} No information on the global load is sent to the users, who know \textit{a priori} that they will pay a fixed price per kWh $p$. Each user consumes energy without any optimization of the system costs and we consider that the consumption profile is given by the original (observed) profile of each user $(\hat{\ell}^h_n)_h$.
The bill of a user $n$ with a total consumption $E_n=\sum_h\sum_a \hat{x}_{na}^h$ will be $b_n^{\text{base}}(\bm{x}) =b_n^{\text{base}}(\bm{x}_n):=p \times E_n$.
As both the PoA \eqref{eq:PoA} and $F$ \eqref{eq: fairness ind} are normalized, the choice of $p$ has no influence at all on the values of those indicators.

\begin{figure}[!t]
\centering
\subfloat{\includegraphics[width=2.7in,height=1.8in,height=1.8in]{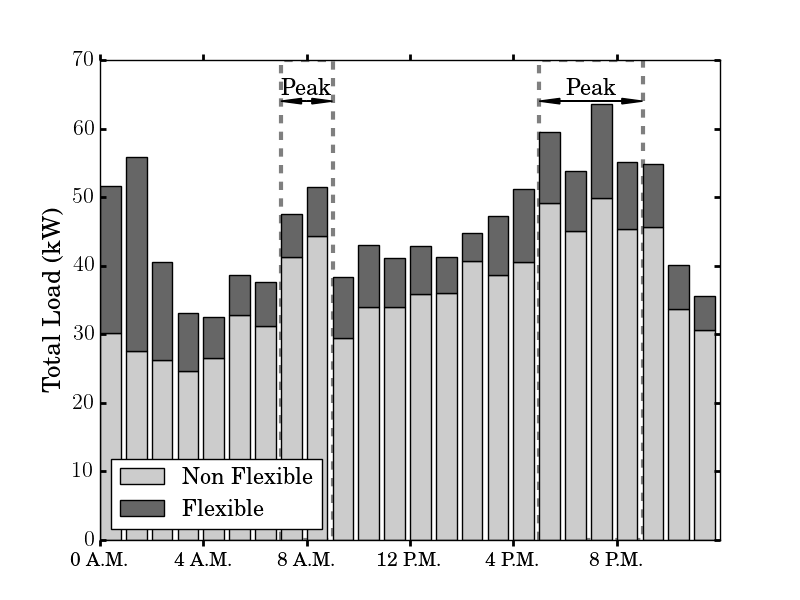}
\label{fig load data}}
\vspace{-0mm}
\subfloat{\includegraphics[width=2.7in,height=1.8in,height=1.8in]{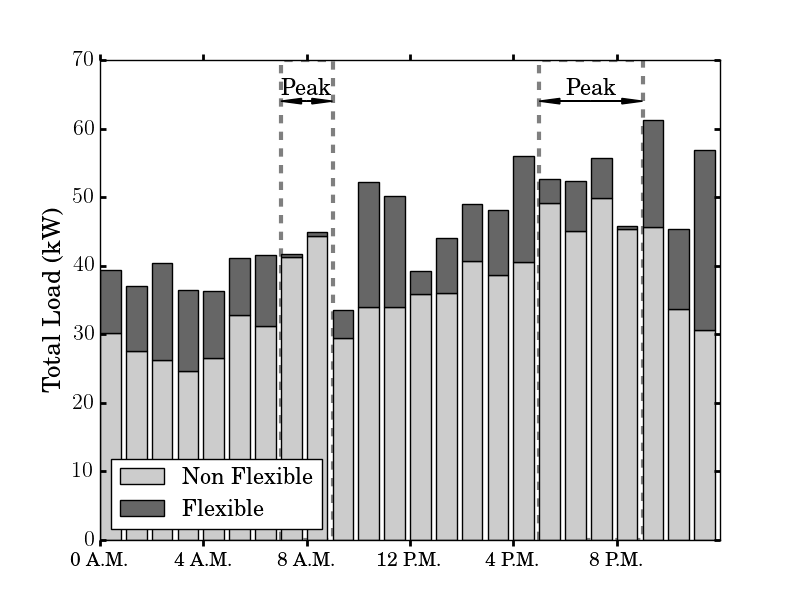}%
\label{fig load peak offpeak}}
\caption{Aggregated load of 30 users on January 10, 2016.\\ \textit{\textbf{top:} Observed profile from the data. \textbf{bottom: } The profile is modified to avoid peak hours in the billing model \eqref{eq: bill peak offpeak}.}}
\label{fig: load data and Peak/offPeak}
\end{figure}

2) \textbf{Peak/Offpeak billing}. This kind of contract already exists in many countries and many of the Texas electricity distributors are proposing  it. The provider defines \textit{a priori} a fixed set of peak hours  $\mathcal{H}_P$ on which the prices are higher. 
We consider that users avoid peak hours as soon as their constraints enable it, by applying a simple greedy algorithm: recursively, a random offpeak hour $h^{\text{off}}$ is chosen and the onpeak load of an appliance $a$ is moved to $h^{\text{off}}$ until $\overline{x}_{na}^{h^{\text{off}}}$ is reached. The resulting consumption profile is denoted by $(\tilde{\ell}^h)_h$ (see \Cref{fig: load data and Peak/offPeak}(bottom) for an example).
 In our simulations, we define $\mathcal{H}_P$ as the set of hours where the nonflexible load is the higher on average, which gives $\mathcal{H}_P= [7\textsc{a.m.}$-$9\textsc{a.m.}]\cup[5\textsc{p.m.}$-$9\textsc{p.m.}] $. 
 We keep the same price ratio $ r^{\text{peak}}= p^{\text{peak}}/p^{\text{off}} =2.84 $ than the Texan distributor \textit{Coserv} \cite{Coserv}. The bill of user $n$ is:
\begin{equation}
\label{eq: bill peak offpeak}
b_n^{\text{Peak/Off}}(\bm{x}_n):=  r^{\text{peak}}p^{\text{off}} \sum_{h\in \mathcal{H}_P} \tilde{\ell}^h_n  +p^{\text{off}} \sum_{h\in \mathcal{H} \setminus \mathcal{H}_P} \tilde{\ell}_n^h \ .
\end{equation}
As explained above, the choice of $p^{\text{off}}$ has no influence at all on the value of the PoA and $F$. However, the ratio $r^{\text{peak}}$ has a direct impact on the fairness indicator \eqref{eq: fairness ind}.

\subsection{Results}

For each day in $\mathcal{D}$, we obtain the equilibrium profile with HP billing \eqref{eq:HP} by running the BRD (\Cref{def: BRD}). In most cases, a hundred iterations (that is, around three optimizations of \eqref{eq:userMinimize}  per user) were sufficient to converge to the equilibrium. The optimal profile (corresponding to the equilibrium in DP) is obtained as the solution of a quadratic program solved with the solver Cplex 12.5. The simulations were implemented in Python 3.5 and run on a single core IntelCore i7-6600U@2.6Ghz with 7.7GB of RAM. The BRD process takes around 50 seconds in average for each simulated day in $\mathcal{D}$,  resulting in a total simulation time of around 20 minutes.

\begin{figure}[!t]
\centering
\subfloat{\includegraphics[width=2.7in,height=1.8in,height=1.8in]{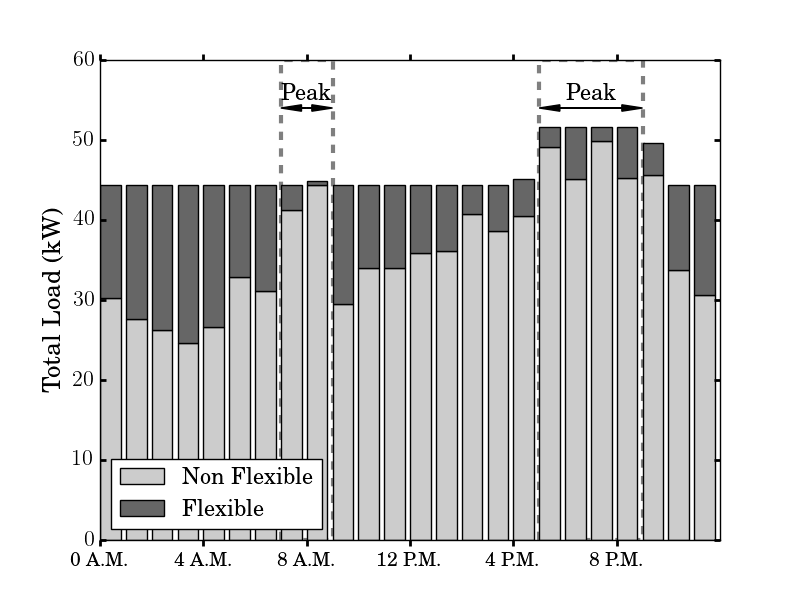}%
\label{fig load opti_first_case}}
\vspace{-0mm}
\subfloat{\includegraphics[width=2.7in,height=1.8in,height=1.8in]{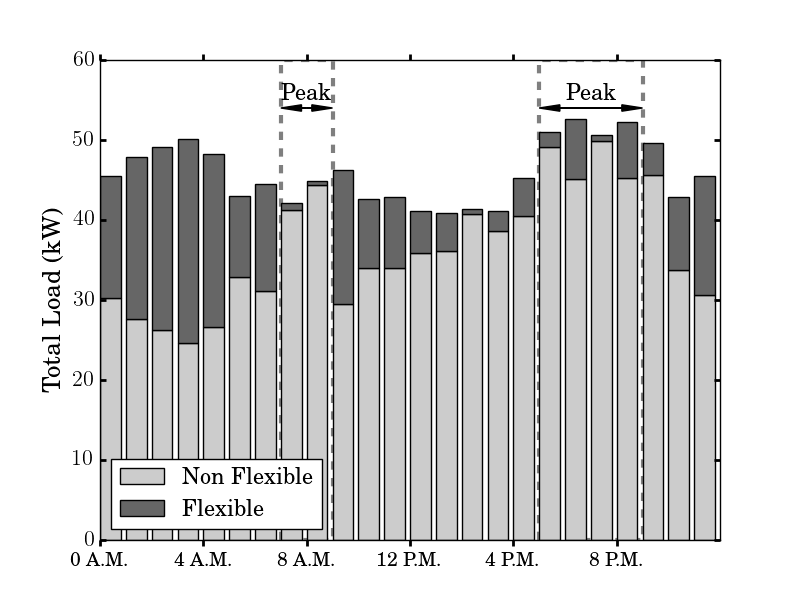}%
\label{fig load opti_second_case}}
\caption{Equilibrium profiles for DP \textit{(top)} and HP\textit{ (bottom)}, on January 10, 2016. \textit{In DP, marginal costs are equal on all hours if the flexible load is sufficient. The HP equilibrium profile remains close to the optimal DP profile.\vspace{-1.5em}}}
\label{fig: optim load HP and DP}
\end{figure}

\begin{figure}[!t]
\centering
\includegraphics[width=3.15in]{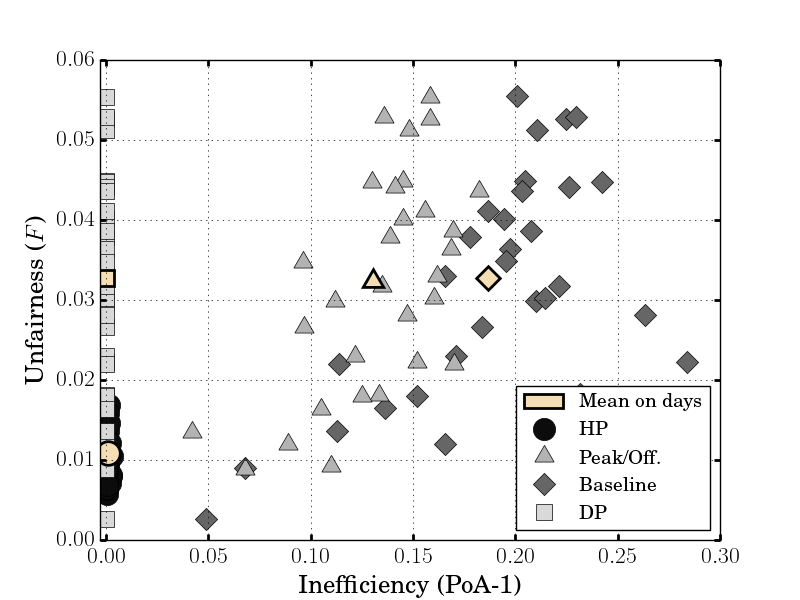}
%
%
\caption{Comparison of billing mechanisms on the 30 days in $\mathcal{D}$ with EV charging considered as flexible. \textit{HP has a PoA of one up to $10^{-3}$. The non-game theoretic billings} Baseline \textit{and} Peak/Offpeak\textit{ are dominated on average. Results are similar when we also consider heating as a flexible appliance.}}
\label{fig fairness eff}
\vspace{-0.5cm}
\end{figure}

\begin{table}
\centering
\resizebox{0.45\textwidth}{!}{
\begin{tabular}{|c||c|c|c|c|}
\hline 

 Flexible Items & \multicolumn{2}{c|}{ EV only}  & \multicolumn{2}{c|}{ EV + furnace} \\
 Billing & PoA-1 (\%) & F (\%)& PoA-1 (\%) & F (\%) \\\hline\hline
HP & 0.0830 (0.0772) & 0.999 (0.286) & 0.0886  (0.104) & 1.17  (0.302)
\\
DP & 0.0  (0.0) & 3.18 (1.38)& 0 (0)	& 3.36  (1.57)\\
Baseline & 18.8 (5.12) & 3.19 (1.38) & 18.5 (5.42) &  3.36 (1.57)
\\
Peak/Off &13.3 (3.17) & 3.20 (1.34) &  12.8 (3.69)	&3.27 (1.04)
\\\hline
\end{tabular}}
\caption{Mean (and standard deviation) of inefficiency (PoA$-1$) and unfairness ($F$) over the days  $\mathcal{D}$ and users $\mathcal{N}$.
}
\label{tab: results}
\vspace{-0.5cm}
\end{table}
The inefficiency (PoA-1) and unfairness ($F$) induced by the four billing mechanisms, that is, DP \eqref{eq: DP}, HP \eqref{eq:HP}, baseline billing and Peak/OffPeak billing \eqref{eq: bill peak offpeak}, are represented in \Cref{fig fairness eff} for each day in $\mathcal{D}$.
The precise values (mean and variance) are given in \Cref{tab: results}. In practice, the PoA of HP is one up to $10^{-3}$: this billing mechanism almost reaches the optimal social cost. The equilibrium profile is not very far from the optimal load profile. \Cref{fig: optim load HP and DP} shows the equilibria of the two mechanisms DP (optimal) an HP. The optimal load profile is very flat. Indeed, due to Kuhn-Tucker conditions of optimality, marginal costs are equal on all hours where constraints \eqref{cons: HP minmax power} are not tight. Therefore, if the part of flexible load is large enough, the total load $L^h$ will be the same for all hours. The equilibrium of the mechanism HP is not as flat, but it remains close to the optimal profile, due to its limited PoA.

Because of the absence of coordination among users in the non-game theoretic Peak/OffPeak billing mechanism, some offpeak hours become congested, as seen on \Cref{fig: load data and Peak/offPeak}(bottom), resulting in high system costs. This efficiency loss is avoided by using a gaming mechanism as HP and DP.

\begin{figure}[!t]
\centering
\includegraphics[width=3.1in]{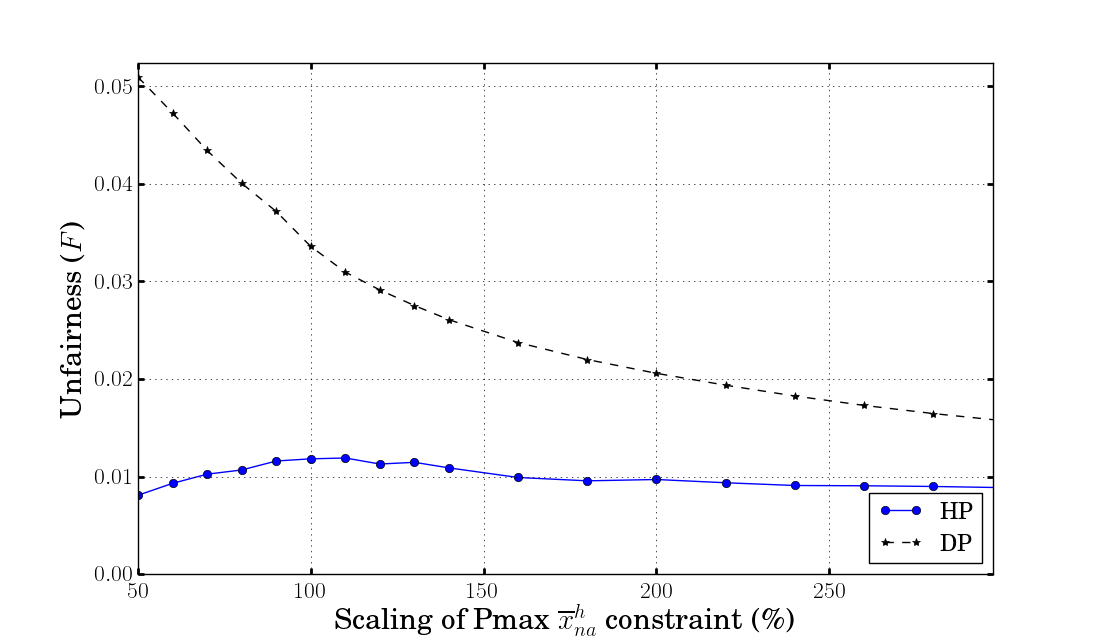}
\caption{Evolution of unfairness in HP and DP with constraints. \\ \textit{When constraints \eqref{cons: HP minmax power} are tight, the DP mechanism has a large unfairness and gets fairer when the constraints  are relaxed.}}
\label{fig: evo unfairness constraint}
\vspace{-0.5cm}
\end{figure}

We can see both from \Cref{fig fairness eff} and \Cref{tab: results} that the HP mechanism achieves an important fairness property in comparison with the other mechanisms. The associated standard deviation of $F$ of 0.3$\%$ indicates that its fairness is also more robust than the other models. Indeed, \Cref{fig: evo unfairness constraint} shows the evolution of the indicator $F$ when we relax the constraint of max power \eqref{cons: HP minmax power} by scaling the value $\overline{x}^h_{na}$ chosen in \Cref{subsec: flex appl} by a factor in $[0.5, 3]$. The unfairness induced by DP decreases when the constraints are relaxed, and it gets closer and closer to HP. Therefore, the HP mechanism will be much more interesting when the constraints are tight.

%


\section{Conclusion}

%
We gave theoretical results ensuring that the hourly proportional billing mechanism has a unique equilibrium and that its price of anarchy is bounded. Experimental results revealed that this mechanism achieves an important fairness property with our quantitative indicator \eqref{eq: fairness ind}, while being very close to the social optimum (up to 0.08\%).
We have seen that the fairness indicator of the hourly mechanism was three times smaller than the other mechanisms, associated with a low variance. As this fairness difference increases with the level of constraints, using the hourly mechanism in practice will be really interesting if the consumption is highly constrained.
If we consider some utility functions in each user optimization (as in \cite{samadi2012advanced}), the daily proportional billing mechanism has no reason to conserve its property of social optimality. Therefore, the small efficiency loss of the hourly mechanism should not have any influence in practice.
%
%


%
%
%



\bibliographystyle{IEEEtran}
\bibliography{../../DRAFTS/Biblio_complete/biblio1,../../DRAFTS/Biblio_complete/biblio2,../../DRAFTS/Biblio_complete/biblioBooks}

\begin{thebibliography}{10}
\providecommand{\url}[1]{#1}
\csname url@samestyle\endcsname
\providecommand{\newblock}{\relax}
\providecommand{\bibinfo}[2]{#2}
\providecommand{\BIBentrySTDinterwordspacing}{\spaceskip=0pt\relax}
\providecommand{\BIBentryALTinterwordstretchfactor}{4}
\providecommand{\BIBentryALTinterwordspacing}{\spaceskip=\fontdimen2\font plus
\BIBentryALTinterwordstretchfactor\fontdimen3\font minus
  \fontdimen4\font\relax}
\providecommand{\BIBforeignlanguage}[2]{{%
\expandafter\ifx\csname l@#1\endcsname\relax
\typeout{** WARNING: IEEEtran.bst: No hyphenation pattern has been}%
\typeout{** loaded for the language `#1'. Using the pattern for}%
\typeout{** the default language instead.}%
\else
\language=\csname l@#1\endcsname
\fi
#2}}
\providecommand{\BIBdecl}{\relax}
\BIBdecl

\bibitem{saad2012game}
W.~Saad, Z.~Han, H.~V. Poor, and T.~Basar, ``Game-theoretic methods for the
  smart grid: An overview of microgrid systems, demand-side management, and
  smart grid communications,'' \emph{IEEE Signal Processing Magazine}, vol.~29,
  no.~5, pp. 86--105, 2012.

\bibitem{siano2014demand}
P.~Siano, ``Demand response and smart grids - a survey,'' \emph{Renewable and
  Sustainable Energy Reviews}, vol.~30, pp. 461--478, 2014.

\bibitem{vardakas2015survey}
J.~S. Vardakas, N.~Zorba, and C.~V. Verikoukis, ``A survey on demand response
  programs in smart grids: Pricing methods and optimization algorithms,''
  \emph{IEEE Comm.s Surveys \& Tutorials}, vol.~17, no.~1, pp.
  152--178, 2015.

\bibitem{mohsenian2010autonomous}
A.-H. Mohsenian-Rad, V.~W. Wong, J.~Jatskevich, R.~Schober, and A.~Leon-Garcia,
  ``Autonomous demand-side management based on game-theoretic energy
  consumption scheduling for the future smart grid,'' \emph{IEEE Trans.
  on Smart Grid}, vol.~1, pp. 320--331, 2010.

\bibitem{chen2010two}
L.~Chen, N.~Li, S.~H. Low, and J.~C. Doyle, ``Two market models for demand
  response in power networks,'' in \emph{Smart Grid Comm.s
  (SmartGridComm), 2010 First IEEE Int. Conf. on}, 2010, pp.
  397--402.

\bibitem{ibars2010distributed}
C.~Ibars, M.~Navarro, and L.~Giupponi, ``Distributed demand management in smart
  grid with a congestion game,'' in \emph{Smart grid communications
  (SmartGridComm), 2010 first IEEE Int. Conf. on}.\hskip 1em plus
  0.5em minus 0.4em\relax IEEE, 2010, pp. 495--500.

\bibitem{li2011optimal}
N.~Li, L.~Chen, and S.~H. Low, ``Optimal demand response based on utility
  maximization in power networks,'' in \emph{Power and Energy Society General
  Meeting, 2011 IEEE}.\hskip 1em plus 0.5em minus 0.4em\relax IEEE, 2011, pp.
  1--8.

\bibitem{baharlouei2012tackling}
Z.~Baharlouei, H.~Narimani, and H.~Mohsenian-Rad, ``Tackling co-existence and
  fairness challenges in autonomous demand side management,'' in \emph{Global
  Comm.s Conf. (GLOBECOM), 2012 IEEE}.\hskip 1em plus 0.5em minus
  0.4em\relax IEEE, 2012, pp. 3159--3164.

\bibitem{fudenberg1991game}
D.~Fudenberg and J.~Tirole, \emph{Game theory, 1991}.\hskip 1em plus 0.5em
  minus 0.4em\relax Cambridge, Massachusetts, 1991.

\bibitem{koutsoupias1999worst}
E.~Koutsoupias and C.~Papadimitriou, ``Worst-case equilibria,'' in \emph{Annual
  Symposium on Theoretical Aspects of Computer Science}.\hskip 1em plus 0.5em
  minus 0.4em\relax Springer, 1999, pp. 404--413.

\bibitem{roughgarden2016twenty}
T.~Roughgarden, \emph{Twenty Lectures on Algorithmic Game Theory}.\hskip 1em
  plus 0.5em minus 0.4em\relax Cambridge University Press, 2016.

\bibitem{johari2005efficiency}
R.~Johari, S.~Mannor, and J.~N. Tsitsiklis, ``Efficiency loss in a network
  resource allocation game: the case of elastic supply,'' \emph{IEEE
  Trans. on Automatic Control}, vol.~50, no.~11, pp. 1712--1724, 2005.

\bibitem{johari2006scalable}
R.~Johari and J.~N. Tsitsiklis, ``A scalable network resource allocation
  mechanism with bounded efficiency loss,'' \emph{IEEE Journal on Selected
  Areas in Comm.s}, vol.~24, no.~5, pp. 992--999, 2006.

\bibitem{christodoulou2005price}
G.~Christodoulou and E.~Koutsoupias, ``The price of anarchy of finite
  congestion games,'' in \emph{Proc. of the thirty-seventh annual ACM
  symposium on Theory of computing}.\hskip 1em plus 0.5em minus 0.4em\relax
  ACM, 2005, pp. 67--73.

\bibitem{roughgarden2006potential}
T.~Roughgarden, ``Potential functions and the inefficiency of equilibria,'' in
  \emph{Proc. of the Int. Congress of Mathematicians (ICM)},
  vol.~3, 2006, pp. 1071--1094.

\bibitem{roughgarden2015intrinsic}
------, ``Intrinsic robustness of the price of anarchy,'' \emph{Journal of the
  ACM (JACM)}, vol.~62, no.~5, p.~32, 2015.

\bibitem{samadi2012advanced}
P.~Samadi, H.~Mohsenian-Rad, R.~Schober, and V.~W. Wong, ``Advanced demand side
  management for the future smart grid using mechanism design,'' \emph{IEEE
  Trans. on Smart Grid}, vol.~3, no.~3, pp. 1170--1180, 2012.

\bibitem{clarke1971multipart}
E.~H. Clarke, ``Multipart pricing of public goods,'' \emph{Public choice},
  vol.~11, no.~1, pp. 17--33, 1971.

\bibitem{moulin2001strategyproof}
H.~Moulin and S.~Shenker, ``Strategyproof sharing of submodular costs: budget
  balance versus efficiency,'' \emph{Economic Theory}, vol.~18, no.~3, pp.
  511--533, 2001.

\bibitem{baharlouei2013achieving}
Z.~Baharlouei, M.~Hashemi, H.~Narimani, and H.~Mohsenian-Rad, ``Achieving
  optimality and fairness in autonomous demand response: Benchmarks and billing
  mechanisms,'' \emph{IEEE Trans. on Smart Grid}, vol.~4, no.~2, pp.
  968--975, 2013.

\bibitem{baharlouei2014efficiency}
Z.~Baharlouei and M.~Hashemi, ``Efficiency-fairness trade-off in
  privacy-preserving autonomous demand side management,'' \emph{IEEE
  Trans. on Smart Grid}, vol.~5, no.~2, pp. 799--808, 2014.

\bibitem{shapley1953value}
L.~S. Shapley, ``A value for n-person games,'' \emph{Contributions to the
  Theory of Games}, vol.~2, no.~28, pp. 307--317, 1953.

\bibitem{rosen1965existence}
J.~B. Rosen, ``Existence and uniqueness of equilibrium points for concave
  n-person games,'' \emph{Econometrica: Journal of the Econometric Society},
  pp. 520--534, 1965.

\bibitem{roughgarden2015local}
T.~Roughgarden and F.~Schoppmann, ``Local smoothness and the price of anarchy
  in splittable congestion games,'' \emph{Journal of Economic Theory}, vol.
  156, pp. 317--342, 2015.

\bibitem{PecanStreet}
\BIBentryALTinterwordspacing
(2016) Pecan street inc. dataport. [Online]. Available:
  \url{https://dataport.pecanstreet.org/data}
\BIBentrySTDinterwordspacing

\bibitem{beaude2016reducing}
O.~Beaude, S.~Lasaulce, M.~Hennebel, and I.~Mohand-Kaci, ``Reducing the impact
  of ev charging operations on the distribution network,'' \emph{IEEE
  Trans. on Smart Grid}, vol.~7, no.~6, pp. 2666--2679, 2016.

\bibitem{Coserv}
\BIBentryALTinterwordspacing
(2017) Coserv electricity distributor. [Online]. Available:
  \url{http://www.coserv.com/Customer-Service/Electric-Rates-And-Tariff}
\BIBentrySTDinterwordspacing

\end{thebibliography}
%

\end{document}